\documentclass[11pt]{article}
\usepackage{latexsym}
\usepackage{amssymb}
\usepackage{amsmath}
\usepackage{amsfonts}
\usepackage{amsopn}
\usepackage{amscd}
\usepackage{graphicx}
\title{Two bijections for sets of words with forbidden factors}
\author{Alexander Valyuzhenich\thanks{Novosibirsk State University, {\tt graphkiper@mail.ru}}}
\begin{document}
\maketitle

\begin{abstract}
In a recent paper by Kitaev and Remmel, several formulas for the number of words of length $n$ avoiding some generalized patterns were established. Each time the obtained function of $n$ had been found in Sloane's Encyclopedia as the number of some other objects, but the bijections between the two sets did not follow from the proof. Kitaev and Remmel stated four open problems on finding respective bijections. Here we solve two of them, concerning sequences A007070 and A048739.
\end{abstract}
Let $\Sigma_k$ denote the alphabet $\{1,2,\ldots,k\}$ and $\Sigma_k^n$ be the set of all words of length $n$ on $\Sigma_k$. In what follows the notation $x=x_1\cdots x_n$ for $x \in \Sigma_k^n$ implies that each $x_i$ is a symbol of $\Sigma_k$. As usual, we say that a word $v$ {\em avoids} a factor $u$ if it cannot be represented as $v=pus$ for any words $p,s$. We also may set other restrictions to words and in particular require the successive letters to differ not too much.

\section{Bijection concerning the sequence A007070}

Let $A_n$ denote the set of all words from $\Sigma_4^n$ avoiding factors 13 and 24, and $B_n$ denote the set of all words $w_1\cdots w_{2n+4}$ from $\Sigma_{7}^{2n+4}$ such that $w_1=1$, $w_{2n+4}=4$, and $|w_i-w_{i+1}|=1$ for all $i=1,\ldots,2n+3$. It was proved in \cite{kit_rem} that $\# A_n=\# B_n$, and both cardinalities are described by the sequence A007070 from \cite{sloan}. However, Problem 2 of that paper was to find a reasonable bijection between the two sets. It is described below.

In fact, we shall build simultaneously the needed bijection between $B_n$ and $A_n$ (denoted by $f$ for all $n$) and the auxiliary bijection $g$ between the set $C_n$ of words $v_1\cdots v_{2n+2}$ from $\Sigma_{7}^{2n+2}$ such that $v_1=1$, $|v_i-v_{i+1}|=1$, and the last letter $v_{2n+2}$ is equal either to 2 or to 6, and the set $D_n$ of words from $\Sigma_4^n$ avoiding factors 13 and 24 and ending with 3 or 4.

For the base of induction, we define both bijections for $n=1,2$. The values of $g$ can be easily listed:
\begin{align*}
g(1212)&=3,\\
g(1232)&=4,\\
g(121212)&=23,\\
g(121232)&=33,\\
g(123212)&=43,\\
g(123232)&=14,\\
g(123456)&=34,\\
g(123432)&=44.
\end{align*}
Now here are the values of $f$ for $n=1$:
\begin{align*}
f(121234)=1,\\
f(123254)=2,\\
f(123234)=3,\\
f(123434)=4.
\end{align*}
Now for $n=2$ there are 14 words in $A_2$, which are all two-letter words except for 13 and 24. Four of them end by 1, four of them end by 2, and six end by 3 or 4.

Now consider the elements of $B_2$: the words of length 8 starting with 1 and ending with 4, with successive letters adjacent. They are also 14. Four of them end by 434: namely, 12345 434, 12343 434, 12123 434, 12323 434; and we somehow bijectively associate with each of them a word from $A_2$ ending with 1. Four others end by 454: namely, 12345 454, 12343 454, 12123 454, 12323 454; they will correspond to elements of $A_2$ ending by 2. The remaining 6 words are 12345654, 12121234, 12123234,
    12321234,
   12343234, 12323234. They somehow correspond to the remaining 6 words from $A_2$. Thus, the values of $f$ on $B_2$ can be defined anyhow under this restriction.

   Now for the induction step let us assume that $f$ and $g$ are already defined on $B_k$ and $C_k$ for all $k<n$ and define them on $B_n$ and $C_n$.

Let us define $g$ first and consider the prefix $x$ of length $2n$ of a word from $C_n$. It could end only by 2, 4, or 6.

If the last letter of $x$ was 4, we define $g(x32)=f(x)23$ and $g(x56)=f(x)14$.

If the last letter of $x$ was 2, we define $g(x12)=g(x)3$ and $g(x32)=g(x)4$.

If the last letter of $x$ was 6, we define $g(x56)=g(x)3$ and $g(x76)=g(x)4$.

It can be easily seen that $g$ is a bijection between $C_n$ and $D_n$ since otherwise either $g$ or $f$ could not be a bijection for shorter words.

Now consider a word $w$ from $B_n$. Clearly, one of the following three situations holds.

Either $w_{2n+2}=4$ and thus $w=w'34$ or $w=w'54$ for some $w'\in B_{n-1}$. Then we define $f(w'34)=f(w')1$ and $f(w'54)=f(w')2$.

Or $w_{2n+2}\neq 4$ but $w_{2n}=4$, so that $w=w''5654$ or $w=w''3234$ for $w''\in B_{n-2}$. Then we define $f(w''3234)=f(w'')23$ and $f(w''5654)=f(w'')14$.

Or $w_{2n+2}\neq 4$ and $w_{2n}\neq 4$, so that the prefix $x$ of $w$ of length $2n$ belongs to $C_{n-1}$. If its last letter was 2, we define $f(x1234)=g(x)3$ and $f(x3234)=g(x)4$. If its last letter was
6, we define $f(x5654)=g(x)3$ and $f(x7654)=g(x)4$.

By the construction, images of distinct words from $B_n$ under $f$ are distinct and belong to $A_n$. Since the cardinalities of the two sets coincide, it is indeed the needed bijection.

\section{Bijection concerning the sequence A048739}

Here the problem is to find the bijection between the set $A_n$ of all words of $\Sigma_3^n$ avoiding factors 13 and $1s3$ for all letters $s$, and the set $B_n$ of all words of $\Sigma_3^{n+3}$ starting with 1, ending with 3, and avoiding factors 13 and 31. The number of elements in any of the sets is described by the sequence A048739 from \cite{sloan}.

Denote by $x_n$, $y_n$, $z_n$ the number of words from $\Sigma_3^{n+3}$ starting with 1, avoiding 13 and 31, and ending by 1, 2, 3, and by $a_n$, $b_n$, $c_n$ the number of elements of $A_n$ ending by 1, 2, 3. We already know from \cite{kit_rem} and \cite{sloan} that $a_n+b_n+c_n=z_n$ and are going to find a bijection $f$ between the respective sets $B_n$ and $A_n$. The base of induction is the following: $f(123)=\lambda$, $f(1123)=1$, $f(1223)=2$, $f(1233)=3$.

Here it is important that $z_n=z_{n-1}+y_{n-1}$ for all $n$ since exactly words ending by 2 or 3 can be extended by 3. Moreover, $a_n=b_n=z_{n-1}$ since any word from $A_n$ can be extended by 1 or 2. Thus, we have $c_n=z_n-a_n-b_n=y_{n-1}-z_{n-1}=x_{n-2}$ since $y_{n-1}=x_{n-2}+y_{n-2}+z_{n-2}$ (any word avoiding 13 and 31 can be extended by 2), and $x_{n-2}=x_{n-3}+y_{n-3}=x_{n-3}+x_{n-4}+y_{n-4}+z_{n-4}=x_{n-3}+x_{n-4}+z_{n-3}$.

So, we also define a bijection $g$ between  the words from $\Sigma_3^{n+1}$ starting with 1, ending with 1 and avoiding 13 and 31, and the elements of $A_n$ ending by 3. For $n=0$, the respective sets are empty, and for $n=1$ we have $g(11)=3$, giving us the base of induction.

Now for the induction step assume that both bijections are already defined for all $k<n$. Define them for $n$ as follows: for $x$ ending by 3 we define $f(x3)=f(x)1$. For $y$ ending by 2 or 3 we define $f(y23)=f(y3)2$ and for $y$ ending  by 1 we define $f(y23)=g(y)$.

It remains to map the $x_{n-2}$ words which are exactly words of $\Sigma_3^{n+3}$ starting with 1, avoiding 13 and 31 and ending by 123 to the remaining $c_n$ words of $A_n$ ending by 3.

If the last letter of $y$ was 1, we define $g(y1)=g(y)3$.

If the last letter of $z$ was 1, we define $g(z21)=g(z)23$.

If the last letter of $z$ was 2 or 3, we define $g(z21)=f(z3)223$.

It can be easily seen now that the bijection $f$ is defined completely and correctly.

\end{document}